\documentclass[12pt,final]{amsart}

\usepackage{showkeys}
\usepackage{xypic}
\newfont{\cyr}{wncyr10}
\newcommand{\Z}{\ensuremath{{\mathbb Z}}}
\newcommand{\Q}{\ensuremath{\mathbb Q}}
\newcommand{\C}{\ensuremath{\mathbb C}}

\newcommand{\coh}[2]{{\rm H}^{#1}(#2)}

\newcommand{\struct}[1]{{\mathcal O}_{#1}}

\newcommand{\spec}{{\rm Spec}}
\newcommand{\dual}{{\hspace{2pt}\check{}}}
\renewcommand{\hom}[2]{{\rm Hom}_{#1}(#2)}

\newtheorem{lem}{Lemma}[section]
\newtheorem{pro}[lem]{Proposition}
\newtheorem{theo}[lem]{Theorem}
\newtheorem{cor}[lem]{Corollary}
\newtheorem{conjecture}[lem]{Conjecture}
\theoremstyle{definition}
\newtheorem{defn}[lem]{Definition}
\newtheorem{rem}[lem]{Remark}
\newtheorem{conv}[lem]{Convention}

\newcommand{\schemes}{{\bf schemes}}
\newcommand{\sets}{{\bf sets}}

\newcommand{\G}{{\mathcal G}}
\newcommand{\bG}{\bar{\mathcal G}}
\renewcommand{\P}{{\mathcal P}}
\renewcommand{\H}{{\mathcal H}}

\newcommand{\dyn}{{\rm Dyn}}
\newcommand{\T}{{\mathfrak T}}
\newcommand{\lie}{{\rm Lie}}
\newcommand{\tw}[2][\G]{{ }^{#2}#1}
\newcommand{\ch}[1][\G]{{\sf X}(#1)}
\newcommand{\chd}[1][\G]{{\sf X}(#1)^\vee}
\newcommand{\gm}{{\mathbb G}_m}
\newcommand{\kk}{{\mathfrak K}}

\newcommand{\Fq}{k}
\newcommand{\Fqn}[1]{k_n}

\newcommand{\Ql}{{\mathbb Q}_{\ell}}
\newcommand{\qql}{{\Ql}}

\newcommand{\tr}{{\rm tr}}
\newcommand{\gal}{{\rm Gal}}

\newcommand{\sha}[1][1]{\text{\cyr X}^{#1}}

\newcommand{\A}{{\mathbb A}}
\renewcommand{\o}{{\mathfrak o}}

\newcommand{\ff}[1][ ]{{K}}
\newcommand{\aff}[1][ ]{\overline{K}}
\newcommand{\cff}[1][ ]{{K_x}}
\newcommand{\acff}[1][ ]{\overline{K_x}}

\newcommand{\gln}[1][n]{{\rm GL}_{#1}}

\newcommand{\bx}{{\bf x}}

\newcommand{\ol}{\overline}

\newcommand{\bun}[2][\G]{{\rm Bun}_{#1}^{#2}}
\newcommand{\sX}{{\mathfrak X}}
\newcommand{\sZ}{{\mathfrak Z}}

\begin{document}
\title[Components of Stacks]{Connected components of moduli stacks of torsors via Tamagawa numbers}

\author{Kai Behrend and  Ajneet Dhillon}  

\subjclass{11E, 11R, 14D, 14H} 

\email{behrend@math.ubc.ca}
\email{adhill3@uwo.ca}

\maketitle

\begin{abstract}
  Let $X$ be a smooth projective geometrically connected curve over
a finite field with function field $K$. Let $\G$ be a connected semisimple group
scheme over $X$. Under certain hypothesis we prove the equality of two numbers associated with $\G$. 
The first is an arithmetic invariant, its Tamagawa number. The second,
is a geometric invariant, the number of connected components of the moduli
stack of $\G$-torsors on $X$. Our results are most useful for studying
connected components as much is known about Tamagawa numbers.
\end{abstract}

\section{Introduction}

We work over a finite ground field $k$. Let $X$ be a smooth geometrically
connected projective curve over $k$ with function field $K$.   Let $\G$ be a semisimple group scheme over $X$.  This means that $\G$ is a smooth group scheme over $X$, all of whose geometric fibres are (connected) semisimple algebraic groups.  We denote the generic fibre of $\G$ by $G$.  

Recall a little bit of terminology: $G$ is {\em split}, if it admits a split maximal torus over $K$.  By the semisimplicity assumption, this implies that $G$ is a Chevalley group, i.e., a group scheme which comes from a split semisimple group defined over $k$ by base extension.    The fundamental group scheme of $\G$ has as fibres the fundamental groups of the fibres of $\G$.  It is a finite abelian group scheme over $X$.  Thus $G$, or equivalently $\G$, is simply connected if and only if this fundamental groups scheme is trivial.

The roots of this article are a circle of ideas that began
to take shape in \cite{harder:70}, \cite{harder:75} and \cite{harder:unpub}.
In particular, it was observed that there is a connection between  Tamagawa numbers
and the trace of the Frobenius endomorphism on the cohomology of certain
moduli spaces. This  led to the following:

\begin{conjecture}[G. Harder] 
If $\G$ is split, the Tamagawa number of $G$ is equal to the number of 
connected components of the moduli space of $\G$-torsors on $X$.
\end{conjecture}

We will study a variation in this paper, that does not assume
that $\G$ is split. We will work 
with stacks rather than spaces. This has two advantages, firstly
there is a precise relationship between the Lefschetz trace 
formula on the moduli stack with the Tamagawa number without the need
for approximations, see \S 4. The second advantage is that these stacks
always exist without any restriction on the characteristic of the
ground field. See for example \cite{balaji:04}.

In order to forge a relationship between
 the stack and the
Tamagawa number we will need to assume that $G$
satisfies the Hasse principle, see \S 4. The Hasse principle is known to hold
for all split groups, see (\ref{C:hassep}). If $\tilde{G}$ is the
universal cover of $G$ then Weil's conjecture asserts that 
the Tamagawa number of $\tilde{G}$ is one. This conjecture is known
in the split case by \cite{harder:74}. In the number field
case the full conjecture is known by \cite{kottwitz:88}.
One expects that a variation of the proof  would work
in the geometric case but this work has not been completed.
The main result of this work is that under certain assumptions
on the ground field the Tamagwa number of $G$ is in fact 
the number of connected components of the moduli stack.
Furthermore one can deduce that these components are in
fact geometrically connected.

In section 3 we begin by recalling the definition of the
Tamagawa number. The section ends with a precise statement
of the main results of this paper.

The purpose of section 4 is to give the proof of the
main results modulo the proof of the trace formula and 
Ono's formula. Section 4.1 recalls basic facts about the
Hasse principle. Section 4.2 use the Seigel formula and the
trace formula to give a geometric interpretation of the 
Tamagawa number. Section 4.3 gives a formula for the
canonical open compact subgroup in terms of special values
of Artin L-functions. Finally the proof is given in 4.4.

Section 5 is devoted to an outline of the proof of the
Lefschetz trace formula for the moduli stack of $\G$-torsors.
The hardest part of the proof is to show that the trace of the 
Frobenius converges absolutely  on the cohomology of the stack.
We only sketch many of the proofs as they are standard (although long) and the
details can be found in \cite{behrend:thesis}. Section 5.1, 
describes the main results on semistabilty for torsors. Section 5.2
introduces the Shatz stratification on the moduli stack of
$\G$-torsors. The proof of the trace formula along with some
semi-purity results for the weight spaces of the cohomology
of the moduli stack of $\G$-torsors are proved in section 5.3.

The final techinical tool needed in section 4.4 is Ono's formula
and some of its consequences. This formula is proved
in section 6.

Let us remark, that Ono's formula implies that the Tamagawa number is equal to the number of elements of $\pi_1(G)$
when $G$ is split.
Thus, we prove that for a semisimple Chevalley group the stack of $\G$-bundles has $|\pi_1(G)|$ components.

In view of \S 4 below, it is tempting to try to use
the results and methods of \cite{atiyah:82} to prove
the main assertions of this paper. However one does
not have the necessary base change theorems required to
transport the results of the cited paper to positive 
characteristic. The moduli stacks of $\G$-torsors are
not proper, indeed they are not even separated. In the case
$\G$ is split and the ground field is $\C$ 
it is known by \cite{atiyah:82} and
\cite{teleman:98} that the moduli stack has
$|{\pi}_1(G,e)|$ components.

\subsection*{Acknowledgements}
This work has its origins in the works of
G\"unter Harder on Tamagawa numbers and moduli spaces of 
bundles. A great intellectual debt is owed to him. 
This paper would not have been possible without
the knowledge of William Casselman. Finally we would
like to thank the referee for for some very useful suggestions
regarding the organisation of this paper.

\section{Notations and Conventions}

\noindent
$\Q_{\ell}$
denotes the $\ell$-adic rationals. We fix once and for all
an inclusion $\Q_{\ell}\hookrightarrow\C$.

\noindent
$k, k_n$  
ground field. Assumed to be perfect. From \S 5 onwards we will 
     assume that $k={\mathbb F}_q$ is a finite field. In this case
$k_n$ denotes the extension of $k$ of degree $n$.

\noindent
$X$
will be a smooth projective geometrically connected curve over $k$ with 
structure map $\pi:X\rightarrow\spec k.$

\noindent 
$\G$ 
is  a smooth connected reductive group scheme over $X$. From \S 6 onwards it will
be assumed to semisimple.

\noindent 
$G$ 
is  the generic fiber of $\G$.. From \S 6 onwards it will
be assumed to semisimple. From (\ref{C:hasse}) onwards we assume
that $G$ satisfies the Hasse principle.

\noindent
$\dyn(\G)$ 
the scheme of Dynkin diagrams of $\G$. See \cite[exp. XXIV]{SGA3}.

\noindent
$\T$ 
the free abelian group on the connected components of $\dyn(\G)$.
 
\noindent
$\T(P)$ 
the free abelian subgroup of $\T$ generated by the components of the type of $P$.
 
\noindent
$\bun[\H,X]{}$ or $\bun[\H]{}$ 
the moduli stack of $\H$-torsors on $X$
where $\H$ is an affine group scheme over $X$.

\noindent
$\bun[\H]{\alpha}$
denotes the moduli stack of $\H$-torsors of degree
$\alpha$.
 
\noindent
$\bun[\H]{\alpha,\le m}$
denotes the moduli stack of $\H$-torsors of degree
of instability at most $m$ and degree $\alpha$.
 
\noindent
$\bun[\H]{\alpha,m}$
denotes the moduli stack of $\H$-torsors of degree
of instability equal to $m$ and degree $\alpha$.
 
\noindent
$\bun[\H]{\alpha,\o}$
denotes the moduli stack of $\H$-torsors of type
of instability $\o$  and degree $\alpha$.
 
\noindent
$\tau(G), \tau_n(G)$
is the Tamagawa number of $G$ or of $G\otimes_{k} k_n$.

\section{The Main Results}

Let us begin by recalling the definition of the Tamagawa number
of a semisimple algebrac group. To do this, we begin by
constructing the Tamagawa measure on $G({\mathbb A})$.

For any point $x$ of $X$ (place of $K$), we denote by $K_x$ the completion of $K$ at $x$. The ring of integers inside $K_x$ is $\widehat{\struct{X,x}}$.  The ring of adeles of $K$, notation $\A$, is the restricted product of all $K_x$ with respect to the $\widehat{\struct{X,x}}$.
Throughout we fix an additive Haar measure $\mu_x$ on each
$K_x$, normalized so that $\widehat{\struct{X,x}}$ has
volume one.  

We fix a section $\omega$ of the line bundle $\wedge_{\struct{X}}^{{\rm dim} G}{\rm Lie}(\G)$.
It induces, in a natural way, 
a Haar measure $\omega_x$ on each of the analytic varieties $\G(K_x)$, 
see \cite[Section 2]{oesterle:84}. The subset $\G(\widehat{\struct{X,x}})$ of 
$\G(K_x)$ is open and its volume is computed by the formula
below, which also characterizes this measure.

\begin{pro}
\label{P:localvol}
  Let $n$ be the order of vanishing of $\omega$ at $x$. Then we
have
\[
{\rm vol}(\G(\widehat{\struct{X,x}})) = |k(x)|^{-n-d}|\G(\Fq(x))|,
\]
where $d$ is the dimension of $G$ and $\Fq(x)$ is the residue field of $x$. 
\end{pro}
 
\begin{proof}
  See \cite[2.5]{oesterle:84}
\end{proof}

For a semisimple group scheme $\G$ on $X$ the vector bundle
$\lie(\G)$ is of degree $0$ on $X$. 

The \textit{Tamagawa measure} is a measure on $G(\A)$ defined by
\[
 q^{(1-g){\rm dim} G}\prod_{x\in X} \omega_x.
\]
It follows from the product formula that this measure does not depend
on the choice of $\omega$. The \textit{Tamagawa number} of $G$, $\tau(G)$,  is
defined to be the volume of $G(\A)/G(K)$ under this measure.
It is known to be finite by \cite{harder:69}.  (Also, the Tamagawa number depends only on the generic fibre $G$ of $\G$, even though we used $\G$ in the definition.)

\begin{conjecture}(Weil)
  If $G$ is simply connected then $\tau(G)=1$.
\end{conjecture}

In the number field case this a theorem proved by
R. Kottwitz, see \cite{kottwitz:88}.

Let $\bun{}$ be the moduli stack of $\G$-torsors on $X$.
We will show that under various hypothesis, the Tamagawa
number computes the number of open and closed substacks
$\bun{}$. Precisely:

\begin{theo}
\label{MT:main1}
Assume that there is a splitting field $L$ for $G$ 
whose constant field is $k$ and that $k$
contains all roots of unity dividing the order
of the fundamental group of $G$. Further assume
that Weil's conjecture holds for the universal
cover of $G$ and $G$ satisfies the Hasse principle.
Then 
$\bun{}$ has $\tau(G)$ components and each of these
components is geometrically connected.
\end{theo}

Note that given any $G$ we can always find $k_n/k$ such
that by base extending to $k_n$ the first two hypothesis of the
theorem are satisfied.

We will deduce from this theorem:

\begin{cor}
\label{MT:cor}
  If the generic fiber of $\G$ is a Chevalley group
then $\bun{}$ has exactly $\tau(G)$ components
each of which are connected.
\end{cor}

Using similar techniques we can also prove:

\begin{theo}
\label{MT:main2}
  If Weil's Tamagawa number conjecture is true then
$\bun{}$ is geometrically connected for every simply connected
$\G$
\end{theo}

\section{The Proof}

\subsection{The Hasse principle}

We begin by recalling some theorems of G. Harder on 
the Hasse principle.  Recall that an algebraic group $G$ over $K$ satisfies the {\em Hasse principle }, if the map of Galois cohomology sets
$$\coh{1}{K,G}\longrightarrow\prod_{x\in X}\coh{1}{K_x,G}$$ is injective.

\begin{theo}[G. Harder]
The Galois cohomology group $\coh{1}{K,G}$ is trivial, if $G$ is simply connected. 
In particular, the Hasse principle holds for such $G$.
\end{theo}

\begin{proof}
  See \cite{harder:75b}.
\end{proof}

\begin{cor}
\label{C:hassep}
 The Hasse principle holds when the generic fiber $G$ is a Chevalley group.
\end{cor}

\begin{proof}
  Let $G'$ be the universal cover of $G$ and $M$ the fundamental group
scheme of $G$. As the first Galois cohomology vanishes for $G'$ and all
its inner forms, we have an injection 
\[
\coh{1}{K,G}\hookrightarrow\coh{2}{K,M}.
\]
As $G'$ is also a Chevalley group, there is an exact sequence
\[
0\rightarrow M\rightarrow \gm^r \rightarrow \gm^r \rightarrow 0,
\]
where $\gm^r$ is a maximal torus of $G'$ containing $M$. The result follows 
from Hilbert's theorem 90 and the fact that there is an injection 
of Brauer groups
\[
{\rm Br}(K)\hookrightarrow\prod_{x\in X}{\rm Br}(K_x).
\]
\end{proof}

\begin{conv}
\label{C:hasse}
  For the remainder of this section we assume that the generic fiber $G$
satisfies the Hasse principle.
\end{conv}

\subsection{The Seigel formula and the trace formula}

We now begin to give a geometric interpretation of the Tamagawa number of $G$.

\begin{lem}
Let $x\in X$ be a closed point.
 The \'{e}tale cohomology set $\coh{1}{\spec(\widehat{\struct{X,x}}),\G_x}$ is trivial.
\end{lem}

\begin{proof}
  We need to show that every $\G_x$-torsor over $\spec(\widehat{\struct{X,x}})$ has a 
$\widehat{\struct{X,x}}$-point. By Lang's theorem such a torsor has a 
point over the residue field of $\widehat{\struct{X,x}}$ which can be lifted to a
$\widehat{\struct{X,x}}$-point by formal smoothness.
\end{proof}

\begin{pro}
\label{P:gtriv}
  Every $\G$-torsor is trivial over the generic point of $X$.
\end{pro}

\begin{proof}
We have a diagram of \'{e}tale cohomology sets
\[
\xymatrix{\coh{1}{X,\G} \ar[r] \ar[d] & \coh{1}{K,G} \ar[d] \\
\prod_{x\in X} \coh{1}{\spec(\widehat{\struct{X,x}}),\G_x} \ar[r] & \prod_{x\in X}
 \coh{1}{K_x,G}.
}
\]
The bottom left corner vanishes and the right vertical map is injective by the
Hasse principle. So the 
top map is trivial.
\end{proof}

Recall that the integral model $\G$ of $G$ defines an open compact subgroup $\kk$ of $\G(\A)$
with 
\[
  \kk = \prod_{x\in X} \G(\widehat{\struct{X,x}}).
\]

\begin{lem}
There is a 
bijection between elements of $G(\A)$ and (isomorphism classes of) triples $(P,\phi,(\rho_x)_{x\in X})$,
where $P$ is a $\G$-torsor, $\phi$ is a trivialization of $P$ over the generic point of 
$P$ and $\rho_x$ is a trivialization of $P$ over the formal disc $\spec(\widehat{\struct{X,x}})$.
\end{lem}

\begin{proof}
  There is an obvious map from such triples to the elements of $G(\A)$. We construct its
inverse as follows. Let ${\bf a}=(a_x)$ be an adelic point of $G$. There is an open
Zariski subset $U$ of $X$ such that ${\bf a}$ is integral over $U$, that is, $a_x\in 
\G(\widehat{\struct{X,x}})$ for $x\in U$.
 Consider the flat cover
\[
U \cup \bigcup_{x\notin U}\spec(\widehat{\struct{X,x}})
\]  
of $X$. To construct $P$ we need only specify descent data with respect to this cover, 
and apply faithfully flat descent. On
the intersection $U\cap \spec(\widehat{\struct{X,x}})$ they are given by $a_x$. 
This gives $P$, together with a generic trivialization and a trivialization at each $x\not\in U$.
For 
$x \in U$ the trivialization $\rho_x$ is given by the trivialization over $U$ multiplied
by $a_x$.
\end{proof}

\begin{pro}
  There is a bijection between points of the double coset space
\[
G(K)\backslash G(\A)/\kk
\]
and the set of isomorphism classes of $\G$-torsors over $X$.
\end{pro}

\begin{proof}
  Use the above Lemma together Proposition (\ref{P:gtriv}).
\end{proof}

\begin{theo}[Siegel's Formula]
\label{T:Siegel}
  We have
\[
\tau(G) = {\rm vol}(\kk) \sum_{P\in {\rm Bun}_{\G}(\Fq)} \frac{1}{|{\rm Aut}(P)|}.
\]
The sum is over isomorphism classes of $\G$-torsors on $X$ and 
$|{\rm Aut}(P)|$ is  the order of the automorphism group of $P$, which
is finite.
\end{theo}

\begin{proof}
We have
\begin{align*}
  \tau(G) & =  {\rm vol}(G(\A)/G(K))        \\
          & =  \sum _x {\rm vol}(\kk xG(K)/G(K))   \\
\intertext{(the sum is over a collection of double coset representatives)}
          & =  \sum _x {\rm vol}(\kk)\frac{1}{|x\kk x^{-1} \cap G(K)|} \\
          & =  {\rm vol}(\kk)\sum _{P\in{\rm Bun}_\G(\Fq)} \frac{1}{{\rm Aut}(P)}.
\end{align*}
The first equality is by the preceeding proposition. One checks in the bijection above 
that the automorphism group of $P$ is identified with $x\kk x^{-1} \cap G(K)$.
Note that one can show that ${\rm vol}(\kk)$ is finite, see \cite{kneser:67}. Furthermore
the sum converges by \cite{harder:69}
\end{proof}

We will prove below a Lefschetz trace formula for the algebraic stack
$\bun{}$.
This formula forges a link between the Siegel formula and the cohomology
of $\bun{}$ and we will describe it now. 

If $\sX$ is an algebraic stack over the finite field $k$ we define
its number of $k$-rational points by
\[
\# \sX(k) = \sum \frac{1}{|Aut(x)|}
\]
where the sum is over isomorphism classes of objects in $\sX(k)$. We denote by 
$\Phi$ acting on the cohomology of $\sX$. We will prove in \S 5 the following
version of the trace formula for $\bun{}$ over $k$ :

\begin{theo}
\label{T:ltfe}
We have 
\[
\sum_{x\in\bun{}(k)}\frac{1}{|{\rm Aut}(x)|} =
q^{(g-1)(\dim G)}\sum (-1)^i{\rm tr}\Phi|_{H^i(\bun{},\qql)}
\]
and both sides converge absolutely.
\end{theo}

\begin{rem}
  The trace formula for stacks of finite type
is proved in \cite{behrend:03}. The stack
$\bun{}$ is not of finite type but it is naturally
filtered by stacks of finite type. Our main
task in proving the above theorem is to prove the convergence.
\end{rem}

\begin{cor}
\label{C:st}
In the current setting 
 we have 
\[
\tau(G) = {\rm vol}(\kk) q^{(g-1)\dim G}\sum_{i=0}^\infty(-1)^i {\rm tr}\Phi|_{\coh{i}{\bun{},\Ql}}. 
\]
\end{cor}

\begin{proof}
  Combine the above with (\ref{T:ltfe}).
\end{proof}

\subsection{Artin $L$-functions and the volume of $\kk$}

We begin by recalling Steinberg's formula for the
number of points of a semisimple group over a finite field.
Let $H/k$ be a semisimple connected linear algebraic group.
By Lang's theorem it is necessarily quasi-split so let $B$
be a Borel subgroup and $T\subseteq B$ a maximal torus.
Let $k_n/k$ be a splitting field for $H$. We form the
Weyl group $W= (N_H(T)/T)(k_n)$. Note that $W$ acts on
$\ch[T\otimes_k k_n]\otimes {\mathbb Q}$ and hence
on its symmetric algebra $S$. By a theorem of Chevalley
the invariants of this action is the symmetric algebra on finite dimensional
graded vector space $V = \oplus_{n\ge 2} V_n$.

\begin{theo}
  Let $F$ be the Frobenius of $k_n/k$. We have
\[
\frac{|H(k)|}{q^d} = \prod_{n\ge 2} \det((1-q^{-n}F)|{V_n}),
\]
where $d$ is dimension of $H$ and $q$ is the number of elements
of $k$.
\end{theo}
 
\begin{proof}
  See \cite[11.16]{steinberg:68}.
\end{proof}

We return to our global situation.  
Firstly observe that:

\begin{pro}
\label{P:volcalc}
  Let $\kk = \prod_{x\in X} \G(\widehat{\struct{X,x}})$ be the
canonical open compact. Then 
\[
{\rm vol}(\kk) = q^{(1-g)\dim G}\prod_x |k(x)|^{-\dim G}|\G(\Fq(x))|.
\] 
\end{pro}

\begin{proof}
  This is by \ref{P:localvol} combined with
the fact that the vector bundle ${\rm Lie}(\G)$ has
degree 0.
\end{proof}

As there is an integral model $\G$ for $G$ we can
by (\ref{l:trick}) find an unramified extension $L/K$ that splits
$G$. We may assume that the extension is 
in fact Galois. Correspondingly we have a Galois cover
\[
Y\rightarrow X.
\]
Again form the
Weyl group $W = (N_G(T)/T)(L)$ which acts on the
symmetric algebra of $\ch[G\otimes_K L]\otimes{\mathbb Q}$. Using
Chevalley's result again we obtain a finite dimensional graded vector
space  $V =\oplus_{n\ge 2} V_n$. Each of the $V_i$'s are
$\gal(L/K)$-modules so we can form the associated
Artin $L$-functions
\[
L_n(X,s) = L(X,V_n,s) = \prod_{x\in X} \det((1 - q^{-s\deg x}f_x)|_{V_n}).
\]
In the above $f_x$ is the Frobenius for the extension of
residue fields $k(y)/k(x)$ and $y$ is an point lying over $x$.
Using the above theorem and (\ref{P:volcalc}) we have
\[
{\rm vol}(\kk) = q^{(1-g)\dim G}\prod_{n\ge 2} L_n(X,n)^{-1}.
\]

In summary:

\begin{theo}
\label{T:volcalc}
The volume of the open compact is given by
\[
{\rm vol}(\kk) = q^{(1-g)\dim G}\prod_{n\ge 2} L_n(X,n)^{-1},
\]
where each $L_n$ is an Artin L-function described above.
Furthermore we have
\[
L_n(X,s) = Z(X,s)^{\gamma_n}\prod_{i=1}^{m_n} p_{ni}(X,s),
\]
where $Z(X,s)$ is the Zeta function and
\[
p_{ni}(X,s) = \prod (1-\alpha_{ni} q^{-s})
\]
and $|\alpha_{ni}|=q^{1/2}$ and $\gamma_n$ is some integer.
\end{theo}

\begin{proof}
  The only part that needs justification is the last
statement. This is by  \cite{milne:80} page 126.
\end{proof}

\subsection{The Proof of the Main Theorems}

We are in a position to give the proof of the main 
theorems modulo some technical results. The proofs of these
will be given later.

 We denote the Tamagawa number of the base
extenstion by
\[
\tau_n(G)=\tau(G\otimes_k k_n).
\]
Here $k_n$ is the unique extension of the finite field
$k$ of degree $n$. A key ingredient in the
proofs of the main theorems  is the fact that the 
sequence
\[
\tau_1(G),\tau_2(G),\ldots,
\]
is constant under suitable hypothesis. Using Ono's formula we will show :

\begin{pro}
\label{P:stable2}
  Suppose that $G$ with field of constants $k$ and
all roots of unity dividing the order of $|\pi_1(G)|$
are in $k$. Further assume Weil's
conjecture for the universal cover of $G$ we have
\[
\tau_n(G) = \tau(G).
\]
for every $n$.
\end{pro}

\begin{proof}
  See (\ref{C:stable2}).
\end{proof}

We need the following technical result to make the proofs
of the main results go more smoothly:

\begin{defn}
  Let $\sum_{m=1}^\infty s_{nm} = t_n$ be a sequence of series of complex
numbers. We say that the series \textit{converge uniformly}
if for every $\epsilon>0$ there is an $M_0$ such that for every
$M\ge M_0$ we have
\[
|\sum_{m=1}^M s_{nm} - t_n| < \epsilon
\]
independently of $n$.
\end{defn}

Our theorem will follow from the following lemma

\begin{lem}
\label{L:silly}
  Let $\sum_{m=1}^\infty s_{nm}$ be a sequence of series that all
sum to $t$ independently of $n$. Furthermore assume that the convergence is uniform
and that the series $\sum_{n=1}^\infty s_{nm}$ converge absolutely for each $m$. Then $t=0$.
\end{lem}

\begin{proof}
  Let $\epsilon>0$ and $M$ be as in the definition of uniform 
convergence. We have
\[
|t - \sum_{m=1}^M s_{nm}| < \epsilon.
\]
However, $\sum_{n=1}^\infty  \sum_{m=1}^M s_{nm}$ converges and hence
\[
\lim_{n\to\infty} \sum_{m=1}^M s_{nm} = 0,
\]
and we are done.
\end{proof}

Recall that the zeta function of $X$ can be written in the form
\begin{eqnarray}
\label{zeta}
Z(X,s) = \frac{\prod_{i=1}^{2g} (1-\alpha_i q^{-s})}{(1-q^{-s})(1-q^{-s+1})},
\end{eqnarray}
where $\alpha_i$ are the eigenvalues of the Frobenius on $\coh{1}{X,\Ql}$. It will
be important below to note that the $\alpha_i$ have absolute value $q^{1/2}$.

Finally, before givig the proof, we will need two facts about the
cohomology of $\bun{}$. Firstly we will need to know that the vector spaces
$\coh{i}{\bun{},\Ql}$ are finite dimensional. Secondly we will need
that the eigenvalues of $\Phi$ have absolute value at most $q^{-i/2}$ on
$\coh{i}{\bun{},\Ql}$. Both these facts are proved in \S5.

\begin{proof}(proof of \ref{MT:main1})
    We have
  \begin{eqnarray*}
    \tau(G)& = &{\rm vol}(\kk)q^{(g-1)\dim G}(\sum_{i=0}^\infty 
     \tr \Phi|_{\coh{i}{\bun{},\Ql}})\quad   (\ref{T:Siegel})\\
           & = & \prod_{n\ge 2}L_n(X,s)^{-1}
      (\tr \Phi|_{\coh{i}{\bun{},\Ql}})\quad   (\ref{T:volcalc}). 
  \end{eqnarray*}
Let $\{\beta_j\}$ be the eigenvalues of $\Phi$ on
\[
\bigoplus_{i>0}\coh{i}{\bun{},\Ql}
\]
and $\epsilon_j$ their signs in the above formula.
So we have
\begin{eqnarray*}
  \tau(G) - \tr \Phi|_{\coh{0}{\bun{},\Ql}} = 
  (\prod_{n\ge 2} L_n(X,n)^{-1} - 1)(\tr \Phi|_{\coh{0}{\bun{},\Ql}}) \\ +
  \prod_{n\ge 2}L_n(X,n)^{-1}(\sum_j \beta_j\epsilon_j).
\end{eqnarray*}
We remind the reader that $\ell$-adic cohomology of a stack
${\mathfrak X}$ over a finite field is defined by first
passing to the algebraic closure, i.e it is really
defined on ${\mathfrak X}\otimes_k {\bar k}$. With this in mind
the action of the Frobenius on the cohomology of the base extension 
${\rm Bun}_{\G_m}$ is just given by $\Phi^m$.
\begin{eqnarray*}
\tau_m(G) - \tr \Phi^m|_{\coh{0}{\bun{},\Ql}} = 
  (\prod_{n\ge 2} L_n(X_m,n)^{-1} - 1)(\tr \Phi^m|_{\coh{0}{\bun{},\Ql}}) \\ 
 +  \prod_{n\ge 2}L_n(X_m,n)^{-1}(\sum_j \beta_j^m\epsilon_j).
\end{eqnarray*}
For future use, we denote the series on the right hand side of the
above equation by $A_m$. 

Note that 
\[
L_n(X_m,s) = L_n(X,s) = Z(X_m,s)^{\gamma_n}\prod_{i=1}^{m_n} p_{ni}(X_m,s)
\]
where the $\gamma_n$ is the same as that in (\ref{T:volcalc}).
We have 
\[
Z(X_m,s) = \frac{\prod_{i=1}^{2g} (1-\alpha_i^m q^{-ms})}
{(1-q^{-ms})(1-q^{m(1-s)})}.
\]
and
\[
p_{ni}(X_m,s) = \prod(1-\alpha_{ni}^m q^{-ms}).
\]
In the above formulas the $\alpha_i$ and $\alpha_{ni}$ 
are the same as those in the formulas for $X$. 
Now $\coh{0}{\bun{},\Ql}$ is finite dimensional by
the results of \S5. So there is a $k_m/k$ such that
the connected compnents of ${\rm Bun}_{\G_m}$ are 
geometrically connected and have a rational point. 
It follows that $\Phi^{lm}$ is the identity on
$\coh{0}{\bun{},\Ql}$ for all $l>1$. So using
(\ref{P:stable2}) the series $A_{lm}$ satisfy 
the first of the hypothesis of (\ref{L:silly}). Now using the
fact that $|\beta_j|\le q^{-1/2}$ by (\ref{C:pure})
the remaining hypothesis are easily checked.
It follows that $A_{lm}=0$. It follows that there
are exactly $\tau(G)$ components. 

Now consider $1<r<m$. A similar analysis shows that each of the
series $A_{lm+r}$ is zero. So $\tr\Phi^{r}|_{\coh{0}{\bun{},\Ql}} = \tau(G)$
also.
It follow that that $\Phi$ must be the
identity and we are done. 
\end{proof}

\begin{proof}(proof of \ref{MT:cor})
  Note that Weil's conjecture is true for Chevalley groups by
\cite{harder:74}. The result is now obtained by
combining the above with (\ref{C:tsplit}).
\end{proof}

\begin{proof}(proof of \ref{MT:main2})
   Argue as in the proof of (\ref{MT:main1}).
\end{proof}

\section{The Lefschetz trace formula for $\bun{}$}

\subsection{Semistability for $\G$-torsors}

The purpose of this subsection is to recall the main results and 
constructions of \cite{behrend:95}. The main point of that paper is
to extend notions such as (semi)stability and Harder-Narasimhan 
filtration to torsors over a reductive group scheme.

For concepts such as root systems with complementary convex solids,
special facets and semistabilty of root systems the reader is referred to
the first three sections of \cite{behrend:95}. The relationship
of these concepts with what is to follow can be found in section 6 of that
paper.

The following construction will be used throughout this work.
\begin{lem}
\label{l:trick}
  There is a finite \'{e}tale cover $f:Y\rightarrow X$ such that
$f^*\G$ is an inner form.
\end{lem}

\begin{proof}
  We make use of the notations of \cite{SGA3}. Let $\G_0$ be 
the constant reductive group scheme over $X$ having the 
same type as $\G$. 
Being an inner form
means that the scheme ${\rm Isomext}(\G,\G_0)$ has a section over $X$.
By \cite[XXIV, theorem 1.3]{SGA3} and  by \cite[XXII, corollary 2.3]{SGA3} $G$
is quasi-isotrivial and hence so is ${\rm Isomext}(\G,\G_0)$. This implies
by \cite[X, corollay 5.4]{SGA3} that ${\rm Isomext}(\G,\G_0)$ is etale 
and finite over $X$. So we take $Y$ to be one of these components
and the section is the tautological section.
\end{proof}

Note that such an inner form is generically split by the Hasse principle. See
(\ref{C:hassep}) and (\ref{P:gtriv}) below.

\begin{defn} 
   Let $\H$ be a smooth affine group scheme over $X$
with connected fibers. We define the \textit{the degree of $\H$} to be
\[
\deg\H = \deg\lie(\H),
\]
where $\lie(\H)$ is the Lie algebra of $\H$ viewed as a vector bundle
on $X$.
\end{defn}

By (\ref{l:trick}) a reductive group scheme has degree 0. 

\begin{defn}
(i) We say that $\G$ is \textit{semistable} if for every parabolic
subgroup $\P$ of $\G$ we have $\deg\P\le 0$.

\noindent
(ii) We say that $\G$ is \textit{stable} if for every parabolic
subgroup $\P$ of $\G$ we have $\deg\P< 0$.

\noindent
(iii) The largest integer $d$, such that there exists a parabolic
subgroup $\P$ of $\G$ of degree $d$ is called the
\textit{degree of instability  of} $\G$ and is denoted
$\deg_i(\G)$.
\end{defn}

By \cite[Lemma 4.3]{behrend:95} the integer $\deg_i(\G)$ is finite.

Let $\dyn(\G)$ be the scheme of Dynkin diagrams of $\G$, see
\cite[XXIV]{SGA3}. The power scheme of $\dyn(\G)$, denoted 
$P(\dyn(\G))$ is the scheme that represents the functor
\begin{eqnarray*}
\schemes/X\rightarrow \sets \\
 T \mapsto {\mathcal P}(\dyn(\G_T)),
\end{eqnarray*}
here ${\mathcal P}$ means set of open and closed subschemes.
For a parabolic subgroup $\P$ of $\G$ recall the definition of
the type of $\P$, denoted $t(\P)$, from \cite[pg. 294]{behrend:95}.
The type $t(\P)$ is a section of $P(\dyn(\G))\rightarrow X$.
In a nutshell,
the type of $\P$ can be thought of in the following way: think of $\G$ as a family
of reductive groups over $X$ and then $\dyn(\G)$ is their Dynkin diagrams glued
together in the appropriate way. Over a point $x\in X$ choose a  Borel subgroup contained inside $\P_x$.
This Borel gives a choice of simple roots which correspond to the vertices of the
Dynkin diagram over $x$. We consider $\lie(\P)_x\subseteq\lie(\G)_x$, and let $R$ be the
subset of the simple roots that consists of those roots $\alpha$ such that the
weight space for $-\alpha$ is in $\lie(P)$ with respect to the above inclusion.
The value of $t(\P)$ over $x$ is the complement of $R$.

Let $\T$ be the free abelian group on $\pi_0(\dyn(\G))$. By
definition of power scheme, the section $t(\P)$ chooses some
connected components of $\dyn(\G)$ and let 
$\T(\P)$ be the free abelian group on these components.

Let $\o$ be a positive element of $\T(\P)$, that is an
element of the form $\sum n_i\o_i$ with the $n_i$ positive.
Given such an $\o$ one can construct a vector bundle 
$W(\P,\o)$, we refer the reader to \cite[pg. 293]{behrend:95}
for the construction and basic properties.

\begin{defn}
  Let $\P$ be a parabolic subgroup of $\G$ and let
$\o$ be a component of its type. We define the \textit{numerical 
invariant of} $\P$ with respect to $\o$ to be 
$\deg W(P,\o)$. The collection of such numbers as $\o$
varies over the components of the type of $\P$ are called  
the \textit{numerical invariants} of $\P$.
\end{defn}

\begin{defn}
  A parabolic subgroup $\P\subseteq\G$ is called
\textit{canonical} if

\noindent
(i) The numerical invariants of $\P$ are all positive.

\noindent
(ii) The Levi component $\P/R_u(\P)$ of $\P$ is
semistable.
\end{defn}

The main results of \cite{behrend:95}
can be summarised in the following.

\begin{theo}
\label{T:main95}
  There is a unique canonical parabolic subgroup of $\G$.
It is maximal among parabolic subgroups of maximal degree.
It commutes with pullback under separable covers.
\end{theo}

The above constructions and definitions apply to a $\G$-torsor $E$
as follows. One forms the inner form
\[
 \tw{E} = E\times_{\G,Ad}\G.
\]
Then $\tw{E}$ is a reductive group scheme over $X$ and we define 
the degree of $E$, etc to be that of $\tw{E}$.

\subsection{The Shatz stratification on $\bun{}$}

We will describe in this subsection the Shatz stratification 
on $\bun{}$ and state its elementary properties. The proofs   are often just generalizations of 
facts about the usual Shatz stratification for vector bundles. 
When new ideas are involved we sketch these. The interested reader is referred to 
\cite{behrend:thesis} for complete proofs.

Let $\bun{}$ be the moduli stack of $\G$-torsors. Let $\ch$
be the group of characters of $\G$. Each $\G$-torsor $E$ defines 
a map
\begin{eqnarray*}
  \deg E & : & \ch \rightarrow \Z \\
         &   & \phi \mapsto \deg(E \times_{\phi}\gm).
\end{eqnarray*}

For $\alpha\in\chd$ we denote by $\bun{\alpha}$ the open and closed
substack of $\bun{}$ of torsors of degree $\alpha$. For $m$ an integer
we denote by $\bun{\alpha,\le m}$ the substack of torsors of degree
of instability at most $m$. It is an open substack of 
$\bun{\alpha}$ that is in fact of finite type. To show
this last fact we proceed in several steps.

By a vector group over $X$ we mean the underlying 
additive group of a vector bundle over $X$.

\begin{pro}
  Let $V$ be a vector group on $X$. Then the natural map
\[
\bun[V]{}\rightarrow\coh{1}{X,V}
\]
makes $\bun[V]{}$ into an affine gerbe over the vector 
space $\coh{1}{X,V}$. This gerbe is trivial, i.e.,
isomorphic to $B\coh{0}{X,V}\times \coh{1}{X,V}$. It follows that
$\bun[V]{}$ is a smooth stack of finite type of dimension
$r(g-1)-d$ where $r$ is the rank of $V$ and $d$ its degree.
\end{pro}

\begin{proof}
  A torsor for $V$ defines, via cocycles, a cohomology class
and this defines the canonical map. It is easy to see
it is a gerbe. let $t:T\rightarrow\coh{1}{X,V}$
be a an affine morphism. The $T$-points of
$\coh{1}{X,V}$ are in bijection with $\coh{1}{X_T, V_T}$,
thinking of $T$ as a $k$-scheme by composition of
structure maps. Hence $t$ defines a cohomology class
$\xi\in\coh{1}{X_T,V_T}$ which corresponds to a $V_T$
torsor $E$. This gives a map $\tilde{t}:T\rightarrow\bun[V]{}$
that lifts $t$. Hence the triviality result.
\end{proof}

\begin{pro}
\label{P:paralevi}
  Let $\P$ be a parabolic subgroup of $\G$ and
let $\H=\P/R_u(\P)$. The natural map
\[
\bun[\P]{}\rightarrow\bun[\H]{}
\]
is a smooth epimorphism of stacks that is of finite
type and relative dimension
\[
\dim_X R_u(\P)(g-1) - \deg(^{E}{\P}),
\]
where $E$ is the universal $\G$-torsor. It induces
an isomorphism on cohomology.
\end{pro}

\begin{proof}
  The unipotent radical is filtered by subgroups
all of whose quotients are vector bundles,
see \cite[XXVI, 2.1]{SGA3}. One then proceeds
by induction.
\end{proof}

\begin{pro}
\label{P:borelfinite}
  Let $B$ be a Borel subgroup of $\G$ and assume that
$\G$ is split over the generic point of $X$. Then
for each $\beta\in \chd[B]$ the stack
$\bun[B]{\beta}$ is of finite type.
\end{pro}

\begin{proof}
  The quotient $B/R_u(B)$ is a split torus.
It is well known that the components of $\bun[\gm]{}$ 
are of finite type and the result follows from the above proposition.
\end{proof}

\begin{pro}
\label{P:cover}
  Let $Z$ be a projective scheme over $k$ and let 
$f:Z'\rightarrow Z$ be a  projective flat cover. Let
$\H$ be a smooth affine group scheme over
$Z$. Then the 
natural pullback map 
\[
\bun[\H]{} \rightarrow\bun[f^*\H]{}
\]
is affine and of finite presentation.
\end{pro}

\begin{proof}
Straightforward.
  See \cite[4.4.3]{behrend:thesis}.
\end{proof}

Before we get to the proof of the fact that 
$\bun{\alpha,\le m}$ is of finite type we 
need some constructions.

Let $\P$ be a parabolic subgroup of $\G$. The type of $\P$
is an open and closed subscheme of $\dyn(\G)$. Its connected
components $\o_1,\o_2,\ldots,\o_s$ generate a subgroup 
$\T(\P)$ of $\T$. There is an action of $\P$ on
$W(P,\o_i)$. Taking the determinant of the 
action  produces a character $\chi_i$ of $\P$.

\begin{defn}
  We say that an element $\alpha$ of $\chd[\P]$ is
positive if $\alpha({\chi_i})>0$ for $i=1,\ldots, s$. 
We denote by $\chd[\P]_+$ the collection of all such 
positive elements. 
\end{defn}

We have a homomorphism 
\[
\T(\P)\rightarrow\ch[\P]
\]
and taking duals and identifying the dual of $\T(\P)$
with itself via the basis $\o_1,\o_2,\ldots,\o_s$ we obtain
\[
\sigma : \chd[\P]\rightarrow\T(\P).
\]
As $P$ acts on $R_u(\P)$ and this group is filtered by 
subgroups with vector bundle quotients we may take 
determinants to obtain a character $\chi_0$. Evaluation at
$\chi_0$ gives a map
\[
m:\chd[\P]\rightarrow\Z.
\] 
Finally the inclusion $\P\subseteq\G$ gives a map
\[
\delta:\chd[\P]\rightarrow\chd.
\]

\begin{pro}
  The map
\[
\delta\times\sigma:\chd[\P]\rightarrow
 \chd\times\T(\P)
\]
is injective with finite cokernel.
\end{pro}  

\begin{proof}
  The details can be found in \cite[7.3.11]{behrend:thesis} but the
idea is as follows. Using (\ref{l:trick}) one can assume
that $\G$ is generically split. To see the reduction observe the cover
in (\ref{l:trick}) may be taken to be Galois with group $\Gamma$. The 
character groups of the original groups are just the groups of $\Gamma$
invariants. 

In the split case one uses the correspondences set up in
\cite[\S 6]{behrend:95} to reduce the question to questions
about root systems with complementary solids.
\end{proof}

Finally:
\begin{theo}
\label{t:finiteness}
  The stack $\bun{\alpha,\le m}$ is of finite type.
\end{theo}

\begin{proof}
See \cite[8.2.6]{behrend:thesis} for full details.
  Again, by passing to Galois covers, 
we may assume that $\G$ is generically split as the natural
map 
\[
\bun{\alpha,\le m}\rightarrow\bun[f^*\G]{tr^\dual(\alpha),\le m}
\]
is of finite type by (\ref{P:cover}). Choose a Borel $B\subseteq\G$
and let $\chi_1,\chi_2,\ldots \chi_s$ be the associated characters.
\[
\coprod_\beta \bun[B]{\beta}\rightarrow \bun[\G]{\alpha,\le m}
\]
where the disjoint union is over all characters such that

\noindent
i. $d(\beta) = \alpha$

\noindent
ii. $m(\delta)\le m$

\noindent
iii. $\beta(\chi_i)\ge -2g $. 

The above Proposition shows this disjoint union is finite. 
A calculation shows that the morphism exists, ie the degree of instability
of the torsors in the image is at most $m$. Furthermore the 
morphism is surjective.  By (\ref{P:borelfinite}) we are 
done.
\end{proof}

Every parabolic subgroup $\P$ of $\G$ determines an element
\[
\sum_{i=1}^s n(P,\o_i)\o_i\quad\text{in}\quad\T,
\]
where $\o_i$ are the connected components of the type of $\P$.
For a reductive group scheme $\G$ on a family of curves 
${\mathfrak X} \rightarrow S$ we define a function 
\[
n:S\rightarrow\T
\] as follows. For a point $s\in S$ choose
an algebraic closure $\overline{k(s)}$ of the residue field
at $s$. Define $n(s)$ to be $n(\P_s)$ where $\P_s$ is the canonical
parabolic subgroup of $\G_s$.

\begin{pro}
  Let $S_d$ be the locally closed subscheme of $S$
where the degree of instability of $\G$ is $d$. Then
$n$ is a continuous function on $S_d$.
\end{pro}

\begin{proof}
  See \cite{behrend:thesis} (7.2.9)
\end{proof}

Denote by $\bun{\alpha,m}$ the locally closed substack of torsors of 
degree $\alpha$ and degree of instability $m$. Denote by 
$\bun{\alpha,\o}$ the locally closed substack of torsors of
degree $\alpha$ and type of instability $\o$.

If $E$ is a $\P$ torsor of degree $\alpha$ then the torsor 
$E\times_\P \G$ has degree $\delta(\alpha)$. If
$\sigma(\alpha) = \sum_{i=1}^s n_i\o_i$ then 
$n_i = n(E\times_{P, Ad} P)$ where we think of
$E\times_{P, Ad} P$ as a parabolic subgroup  of
$E\times_{P, Ad} G$. Furthermore, $\deg_i(E\times_{P, Ad} P)) = m(\alpha)$.

\begin{theo}
\label{T:radical}
  Denote by $\bG$ the reductive group scheme $\G\times_k \bar{k}$ over
the curve $X\times_k \bar{k}$. Let $\P$ be a parabolic subgroup of
$\bG$ and let $\alpha\in\chd[\P]_+$.
  The natural map
\[
\bun[\P]{\alpha, 0} \rightarrow \bun[\bG]{\delta(\alpha), m(\alpha)}
\]
is finite radical and surjective.
\end{theo}

\begin{proof}
Recall that a morphism is radical if induces a bijection
on $L$-points for every field $L$.
  Representability of this morphism is easy to show.
  The fact that it is radical and surjective  amounts to the 
existence and  uniqueness  of the canonical 
parabolic. 
\end{proof}

\subsection{The Lefschetz trace formula for $\bun{}$}

\begin{pro} 
\label{P:iso1}
  Let $\P$ be a parabolic subgroup of $\bG$. Let $\alpha\in \chd[\P]_{+}.$
Let $\H = \P/R_u(\P)$.
Then there is a natural isomorphism 
\[
\coh{i}{\bun[\bG]{d(\alpha), \sigma(\alpha)},\Ql} \rightarrow
\coh{i}{\bun[\H]{\alpha,0}}.
\]
\end{pro}

\begin{proof}
  Use (\ref{T:radical}) and (\ref{P:paralevi}). Note
that a finite radical and surjective morphism induces
an isomorphism on cohomology. For this last fact
see \cite[Expose VII]{SGA4}
\end{proof}

\begin{lem}
\label{L:r}
  There is a 
function $r:\T\rightarrow\Z$ such that if $E$ 
is a $\bG$-torsor of type of instability $\o$ then 
$r(\o) = \dim_X R_u(\P)$ where $\P$ is the canonical
parabolic of $E$.
\end{lem}

\begin{proof}
  This is just because two parabolic subgroups having the
same type are twisted forms of each other.
\end{proof}

\begin{pro}
\label{P:iso2}
  The closed immersion
\[
\bun{\alpha,\o}\rightarrow\bun{\alpha,\le m(\o)}
\]
is of codimension $c(\o)=r(\o)(g-1)+m(\o)$.
\end{pro}

\begin{proof}
  This is a standard dimension calculation.
\end{proof}
 
Define $\gamma(i)$ to be the smallest integer such that
\[
\gamma(i) \ge \left\{
  \begin{array}{rr}
     1 + i/2 & \text{ if } g>0 \\
     1 + i/2 + |\Phi| & \text{ if } g=0
  \end{array}\right.
\]
where $|\Phi|$ is the number of roots of $\G$.

\begin{pro}
\label{P:iso3}
Let $i\ge 0$ be such that $m\ge \gamma(i)$. Then the
canonical map
\[
\coh{i}{\bun{\alpha,\le m},\Ql}\rightarrow
\coh{i}{\bun{\alpha},\Ql}
\] 
is an isomorphism.
\end{pro}

\begin{proof}
 Let $c$ be the codimension of $\bun{\alpha,\o}$
in $\bun{\alpha,\le m}$ where $m(\o)=m$.
  Using the above one shows that for 
$m\ge \gamma(i)$ we have $i\le 2c-2$. A Gysin
sequence yields the result.
\end{proof}

For a divisor $D$ on $X$ we denote by 
$\bun{}(D)$ the moduli stack of $\G$-torsors
with level structure over $D$. By \textit{level structure}
over $D$ for a torsor $E$ we mean a section of $i_D^*E$ where
\[
i_D:D\hookrightarrow X
\]
is the natural inclusion. 

\begin{pro}
\label{P:DM}
  Let $D$ be a divisor on $X$.  Let $\bun{\alpha,\le m}(D)$ be the
moduli stack of $G$-torsors on $X$ with level structure at $D$. Let $E$ 
be the universal torsor on $\bun{\alpha}$. If 
\[
p_* \tw[\lie(\G)]{E}(-D)|_{\bun{\alpha,\le m}} =0 
\]
then $\bun{\alpha,\le m}(D)$ is a Deligne-Mumford stack.
\end{pro} 

\begin{proof}
 We need to show that the diagonal morphism  of 
$\bun{\alpha,\le m}(D)$ is unramified.
  Let $K$ be a field and consider a $K$-point
\[
\spec K \rightarrow \bun{\alpha,\le m}(D).
\]
It corresponds to a pair $(E,s)$ where $E$ is a
a torsor on $X_K$ and $s$ is a section. Denote by
${\rm Aut}(E,s)$ the automorphism group of $E$ compatible
with $s$. If $\pi:X\rightarrow\spec k$ is the structure
morphism then we need to show that 
$\pi_K({\rm Aut}(E,s))$ is unramified over $\spec K$.
A calculation shows that one may identify the
Zariski tangent spaces of this group with
$\coh{0}{X_K, \lie(\G)(-D)}$ which vanishes. 
\end{proof}

\begin{theo} 
\label{T:pure}
  The eigenvalues of the arithmetic Frobenius acting on 
$\coh{i}{\bun{\alpha,\le m},\Ql}$ have absolute value at most
$q^{-i/2}$.
\end{theo}

\begin{proof}
  First, for a smooth Deligne-Mumford stack 
the analogous statement is true by comparison to 
its coarse moduli space. See \cite{behrend:93}.
Now for a quotient stack $[\sX/\Gamma]$
where $\sX$ is a smooth Deligne-Mumford stack
and $\Gamma$ any algebraic group, we proceed as follows. First 
we may assume that $\Gamma=\gln$ by choosing
a faithful representation 
\[
\Gamma\hookrightarrow\gln
\]
and replace $\sX$ by $\sX\times_{\Gamma}\gln$.
The Leray spectral sequence for 
\[
[\sX/\gln]\rightarrow B\gln
\]
has $E_2$ term
\[
\coh{i}{B\gln,\Ql}\otimes\coh{j}{\sX,\Ql}.
\]
The result follows for the quotient as it is 
known for $\sX$ and $B\gln$.
\end{proof}

\begin{cor}
 \label{C:pure}
  The eigenvalues of the arithmetic Frobenius acting on
$\coh{i}{\bun{\alpha}}$ have absolute value at most
$q^{-i/2}$.
\end{cor}

\begin{lem}
  For the group $\bar{\G}$ there is a function
\[
m:\T\rightarrow\Z
\]
such that if $E$ is a $\bar{\G}$-torsor with type of
instability $\o\in\T$ then $m(\o)$ is the degree
of instability of $\tw[\bar{\G}]{E}$. 
If $\P$ is a parabolic subgroup of $\bar{G}$ then the
diagram
\[
\xymatrix{
\chd[\P]  \ar[r]^-\o  \ar[dr] & \T \ar[d]^m \\
 & \Z}
\]
commutes.
\end{lem}

\begin{proof}
  If $E$ and $F$ are torsors with the
same type of instability $\o$. There  is a parabolic
$\P$ of type $\eta$ where $\eta$ is the support of $\o$ 
and $E$ and $F$ have reductions $E'$ and $F'$ to $\P$. Then
$\o(\deg E') = \o = \o(\deg F')$ and our lemma follows from
the fact that $m$ factors through $\T(\P)$. 
\end{proof}

Let $\eta$ be a closed and open subscheme of $\dyn(\bar{\G})$. 
Let 
\[
C(\eta,mu) =|\{ \o\in \T(\eta)^+ | m(\o) = \mu \}|.
\]
Here $\T(\eta)^+$ denotes the set of linear combinations in the support of
$\eta$ all of whose coefficients are positive.

\begin{lem}
  We have $C(\eta,\mu) = O(\mu^s)$ where $s$ is the number of
components of the type of $\P$.
\end{lem}

\begin{proof}
  Recall the definition of the characters $\chi_i$ and
$\chi_0$ from the discussion after Proposition (\ref{P:cover}).
The result follows from the fact that there are positive rational numbers such that
\[
\chi_0 = \sum_{i=1}^s y_i \chi_i.
\]
\end{proof}

\begin{lem}
  Let $R$ be the radical of $\bar{\G}$. There exist
finitely many $d_1,\ldots, d_n\in \chd[\bar{\G}]$ such that
for every $d\in\chd[\bar{\G}]$ there is an $R$-torsor of
degree $d_i-d$.  
\end{lem}

\begin{proof}
  Let 
\[
M = \{ \delta\in\chd[R] | \bun[R]{\delta}\ne\emptyset \}.
\]
By looking at dual objects one constructs an exact sequence of
tori over $X$
\[
1\rightarrow S\rightarrow R\rightarrow\gm^r \rightarrow 1
\]
where the last map has a quasi-section. We can identify 
$ \chd[R]$ with $\chd[\gm^r]$ and using this identification and 
the quasi-section we observe that $M$ has finite index in 
$\chd[R]$. It follows that $M$ has finite index in 
$\chd[\bar{\G}]$ and we take $d_i$ to be a set of
coset representatives for $\chd[\bar{\G}]/M$.
\end{proof}

Let $A$ be the set of open and closed subschemes 
$\eta\subseteq\dyn\bar{\G}$ such that there is a 
torsor $E$ whose canonical parabolic has type $\eta$.
For each $\eta\in A$ fix a torsor $E_\eta$. Let
$P_\eta\subset\tw[\bar{\G}]{E_\eta}$ be the canonical 
parabolic and let $H_\eta=P_\eta/R_u(P_\eta)$ be its
Levi factor. We choose for each $\eta$ degrees 
$d(\eta,1) , d(\eta,2),\ldots, d(\eta,n)\in\chd[H_\eta]$
according to the above lemma. We get a finite family
of stacks
\[
\bun[H_\eta]{d(\eta,j),0}
\]
parameterized by $A\times\{1,\ldots,n\}$. Set
\[
b_i(\eta,j) = \dim_{\Ql}\coh{i}{\bun[H_\eta]{d(\eta,j),0},\Ql}.
\]
Choose $B(i)= \sup b_i(\eta,j)$. 

\begin{lem}
  There is an integer $N$ so that $B(i)=O(i^N)$.
\end{lem}

\begin{proof}
  The stacks in question are quotients of smooth 
Deligne-Mumford stacks by (\ref{P:DM}). The result follows
from the usual spectral sequence.
\end{proof}

Recall the definitions of $\gamma(i)$  after (\ref{P:iso2}) 
and $r(\o)$ in (\ref{L:r}). We set
\[
D(i) = \sum_{eta\in A} \sum_{\mu=0}^{\gamma(i)} C(\eta,\mu) B(i-2\mu-2r(\eta)(g-1)).
\]

\begin{lem}
The following sum converges:
\[
\sum_{\eta\in A}\sum_m q^{-m + (1-g)R_uP} \sum_i \dim 
H^i(\bun{H_\eta,C}{d(m),0},\Ql) q^{-i/2}.
\]
\end{lem}

\begin{proof}
One use the above estimate.
\end{proof}

The convergence now follows from some general observations about
cohomology of stacks that we outline below. Denote by
$W^iH^j(X,\Ql)$ the $i$ th weight space of the $j$th cohomology
group. 

Recall that a morphism $Z\rightarrow \tilde{Z}$ is called
a universal homeomorphism if it is finite radical and surjective.
Such a morphism induces an equivalence of etale sites by
\cite[Expose IX 4.10]{SGA1}. 

\begin{theo}
Let $\sX$ be a smooth stack with a countable stratification 
by locally closed stacks $\widetilde{\sZ}_i$. We assume that the union
$\sX_n=\cup_{i=0}^n\widetilde{\sZ_i}$ is open. 
Suppose that there
are universal homeomorphisms $\sZ_i\rightarrow\widetilde{\sZ}_i$ 
with each of the $\sZ_i$'s smooth and the following sum
converges:
\[
\sum_{n=0}^\infty q^{-{\rm codim}(\sZ_n,\sX)}\sum_{i,j} 
\dim Gr^W_iH^j(\ol{\sZ}_n,\qql)q^{-i/2}<\infty\,.
\] 
Then the trace of the Frobenius converges absolutely
on $\sX$.
\end{theo}

\begin{proof}
Let  $Z\to X$ be a morphism of finite type smooth
schemes which factors as $ Z\to\widetilde Z\to X$, where
$\pi:Z\to\widetilde Z$ is a universal homeomorphism and
$i:\widetilde Z\to X$ a closed immersion with complement $U$. We
have a long exact sequence
$$\ldots\to H^\ast(\widetilde Z,i^!\qql)\to 
H^\ast( X,\qql)\to 
H^\ast( U,\qql)\to\ldots
$$
Let $c=\dim X-\dim Z$.  We have
$$H^{\ast-2c}(Z,\qql(-c))=H^\ast(Z,\pi^!i^!\qql)\,$$
because $Z$ and $X$ are smooth.  Now pulling back via $\pi$ induces an
isomorphism of \'etale sites (see \cite[Expose IX,4.10]{SGA1}).
As $\pi_\ast$ is  the
right adjoint of $\pi^\ast$,  it is the inverse of $\pi^\ast$ and hence also
a left adjoint of $\pi^\ast$.  Since $\pi$ is proper, 
we conclude that $\pi^!=\pi^\ast$. Thus, we have
$$H^\ast(Z,\pi^!i^!\qql)=H^\ast(Z,\pi^\ast i^!\qql)=H^\ast(\widetilde
Z,i^!\qql)\,,$$
Thus we have a natural long exact sequence
$$\ldots\to H^{\ast-2c}( Z,\qql(-c))\to 
H^\ast( X,\qql)\to 
H^\ast( U,\qql)\to\ldots
$$
This result extends to stacks and filtrations of schemes and stacks
consisting of more than two pieces.
Assembling these long exact sequences we get the required result
and some simple analysis gives the required result.
\end{proof}

\begin{theo}
\label{T:ltf}
We have 
\[
\sum_{x\in\bun{\alpha}(k)}\frac{1}{{\rm Aut}(x)} =
q^{(g-1)(\dim G)}\sum (-1)^i{\rm tr}\Phi|_{H^i(\bun{},\qql)}
\]
and both sides converge absolutely.
\end{theo}

\begin{proof}
The convergence of the trace is by the aboove proposition and
lemma, noting that the natural maps 
\[
\bun[P]{m,0}\rightarrow\bun{}
\]
are finite radical and surjective onto their image
by the uniqueness of the canonical parabolic. The stratification
being used is Shatz stratification induced by reduction to
the canonical parabolic.
The open substacks
of bounded degree of instability are of finite type by 
\ref{t:finiteness}. As the trace formula holds for these
the result follows.
\end{proof}

\section{Ono's Formula and Applications}

Let $M$ be the fundamental group scheme of $G$. So $M$ is an Abelian group scheme
over $K$ and there is an exact sequence
\[
1\rightarrow M \rightarrow \tilde{G} \rightarrow G \rightarrow 1,
\]
where $\tilde{G}$ is the universal cover of $G$. For a
continuous $\gal(\aff^s/\ff)$-module $N$ 
let
\[
\sha(K,N) = \ker\big(\coh{1}{\ff,N}\rightarrow
\prod_{x\in X}\coh{1}{\cff),N} \big)
\]
Here $\cff$ is the completion of the global field $K$ 
at $x$.
\begin{theo}(Ono's Formula)
\label{T:ono}
Assume that Weil's conjecture holds true for the 
universal cover $\tilde{G}$ of $G$, that is
\[
\tau(\tilde{G})=1.
\]
Then we have
\[
\tau(G) = \frac{|\coh{0}{\ff,\widehat{M}}|}
{|\sha(\widehat{M})|},
\]
where $\widehat{M}=\hom{}{M,\gm}$ is the dual Galois module.
\end{theo}

A notational clarification is in order here.
The object  $\widehat{M}$ is to be viewed as functor
on field extensions of $K$. We have
\[
\hat{M}(L) = \ch[M] \otimes_K L.
\]
We will often write $\coh{i}{L,\hat{M}}$ when we really mean
$\coh{i}{L,\hat{M}(\bar{L}^{sep})}$.

The above theorem is the main result of \cite{ono:65}. It was originally only proved
in the number field case and some modifications are needed in the
function field case. We will detail these below.

To prove this result we need to generalize the theory of Tamagawa measures
to reductive groups. We refer the reader to
 \cite[1.4]{oesterle:84} for 
the definition of the Tamagawa measure $d\tau_H$ for a reductive
group $H$.

In \cite{ono:65}, the theorem is proved by reducing to the case
of an isogeny of tori. This was treated in \cite{ono:63}. However
this last paper contains a small error in the function field case
which was corrected in \cite[pg.23 and Chapter IV]{oesterle:84}.

We need some background results  before giving the proof 
of the above theorem.
Let $\Lambda_1\subseteq\Lambda$ be an inclusion of free abelian
groups of the same rank $r$. Let 
\[
\bx = \{x_1 =0, x_1,\ldots x_t\}
\]
be a set of coset representatives for $\Lambda/\Lambda_1$. A function 
$f:\Lambda\rightarrow {\mathbb R}$ is said to be $\bx$-\textit{compatible} 
if 

\noindent
(i) $f$ has finite support.

\noindent
(ii) $f(\alpha) = f(\alpha + x_i)$ for all $\alpha\in \Lambda_1$ and all $i$.

Such a function is said to be $\Lambda_1$-\textit{compatible} if 
it is $\bx$-compatible for some choice of coset representatives containing $0$.

\begin{lem}
\label{L:easy}
  Suppose we have three lattices 
\[
\Lambda_1\subseteq\Lambda_2\subseteq\Lambda_3
\]
of the same rank $r$. Let $f:\Lambda_3\rightarrow{\mathbb R}$
be $\Lambda_1$ compatible then
\[
\sum_{y\in \Lambda_2} f(y) = (\sum_{y\in \Lambda_3}f(Y)) \frac{1}{[\Lambda_3:\Lambda_2]}.
\]
\end{lem}

\begin{proof}
  This is elementary.
\end{proof}

We denote by $D_q$ the set $\{ q^{i}|i\in \Z\}$. Choose
a basis $\{\chi_1,\ldots\chi_r\} = \chi$ for the group $\ch[H]$ of rational
characters of $H$. We define
\[
\psi_H^{\bf \chi} = \psi_H : H(\A) \rightarrow D_q^r
\]
by sending 
\[
x\mapsto (\| \chi_1(x)\|,\| \chi_2(x)\|,\ldots, \|\chi_r(x)\|).  
\]
In the above $\chi_i$ is really the adelization of $\chi_i$. 
The image of $\psi$ is of finite index in $D_q^r$. Denote
by $H(\A)^1$ the kernel of $\psi$. We denote by $d\tau^1_H$ the
measure on $H(\A)^1$ that is the quotient of $d\tau_H$ by
$[D_q^r:Im(\psi)](\log q)^r$. The Tamagawa number of $H$ is
\[
\tau(H)=\int_{H(\A)^1/H(K)} d\tau^1_H.
\]

\begin{pro}
\label{P:compatible}
Let $\Lambda\subseteq Im(\psi) \subseteq D_q^r$ be a 
sublattice of maximal rank. If $F$ is $\Lambda$-compatible then
\[
\int_{H(\A)/H(K)} F(\psi_H(x))d\tau_H = \tau(H) 
(\sum_{x\in D_q^r} F(x))(\log q)^r
\]  
\end{pro}
\begin{proof}
  We have, noting the previous lemma,
\begin{eqnarray*}
 \int_{H(\A)/H(K)} F(\psi_H(x))d\tau_H &=& \sum_{y\in Im(\psi)} F(y)\int_{H(\A)^1/H(K)} d\tau \\
       &=&  \sum_{y\in D_q^r} F(y)\int_{H(\A)^1/H(K)} \frac{d\tau}{[D_q^r:Im(\psi)]} \\
       &=&  \sum_{y\in D_q^r} F(y)\int_{H(\A)^1/H(K)} {d\tau^1_H}(\log q)^r.
\end{eqnarray*}
\end{proof}

\begin{pro}
\label{P:compare1}
  Consider an exact sequence
\[
1\rightarrow H' \stackrel{i}{\rightarrow}H 
\rightarrow H'' \rightarrow 1.
\]
Suppose $H'$ is a torus, $H''$ is semisimple and the sequence is 
generically split. Then 
\[
\tau(H')\tau(H'') = \tau(H)|{\rm cok}\hat{i}|.
\]
In the above formula $\hat{i}$ is the dual map on character
groups.
\end{pro}

\begin{proof}(cf. \cite[Proposition(1.2.2)]{ono:65})
  The fact that the sequence is generically split implies   that 
the induced map on adelic and $K$-points is exact. Furthermore
$\ch[H]$ is a subgroup of $\ch[H']$ of finite index $|{\rm cok}\hat{i}|$.
By the elementary divisors theorem we may choose a basis 
$\chi_1 \ldots \chi_r$ of $\ch[H']$ such that 
 $m_1\chi_1 \ldots m_r\chi_r$  is a basis of $\ch[H]$.
We have a diagram
\[
\xymatrix{
D_q^r \ar[r]^{\textbf{m}} & D_q^r \\
Im(\psi_{H'}) \ar[r] \ar[u] & Im(\psi_H) \ar[u]\\
H'(\A) \ar[u]_{\psi_{H'}} \ar[r] & H(\A) \ar[u]^{\psi_H}.
}
\]
We have a sequence of inclusions
\[
Im(\psi_{H'}) \hookrightarrow D_q^r \hookrightarrow D_q^r.
\]
In what follows we write the group operation on $D_q^r$ additively.
There are basis of the two outside lattices of the form
\[
\{e_1,\ldots e_r \}\quad \{ d_1e_1,\ldots d_re_r \},
\]
by the elementary divisors theorem.
Define a function on $D_q^r$ by 
\[
f(e_1^{\alpha_1} + \ldots e_r^{\alpha_r}) = 
\left\{
  \begin{array}{cc}
      1 & 0\le \alpha_i < d_i \\
      0 & \text{otherwise}
  \end{array}\right.
\]
This function has the  property that
\[
\int_{H'(\A)/H'(K)} f(\psi_H(x') + t) d\tau_{H'} = \int_{H'(\A)/H'(K)} f(\psi_H(x')) d\tau_{H'},
\]
which just follows from the definitions.
Using the compatibility properties of $f$ and (\ref{L:easy}) one shows that
\[
\int_{H'(\A)/H'(K)} f(\psi_H(x')) d\tau_{H'} = \frac{1}{|{\rm cok}\hat{i}|}
 \int_{H'(\A)/H'(K)} f(\psi_{H'}(x')).
\]
This follows from the definition of $f$ and the fact that $d\tau_{H'}$ is a Haar measure. 
We have 
\begin{eqnarray*} 
  \tau(H)(\sum_{y\in D_q^r} f(y))(\log q)^r &=&  \int_{H(\A)/H(K)} f(\psi_H(x))d\tau_H \\
        & =& \int_{{H''}(\A)/{H''}(K)} d\tau_{{H''}}\int_{H'(\A)/H'(K)} f(\psi_H(x')\psi(x)) d\tau_{H'}\\
        &=& \int \tau({H''})\int_{H'(\A)/H'(K)} f(\psi_H(x')) d\tau_{H'} \\
        &=&  \frac{\tau(H')\tau({H'}')}{|{\rm cok}\hat{i}|}(\sum_{y\in D_q^r} f(y))(\log q)^r,
\end{eqnarray*} 
which finishes the proof.
\end{proof}

\begin{lem}
  Let 
\[
1\rightarrow H' \rightarrow H 
\stackrel{\kappa}{\rightarrow} 
{H''} \rightarrow 1
\]
be an exact sequence of linear algebraic groups over $K$. Assume
that $H'$ is semisimple simply connected, ${H''}$ is a torus and $H$ is reductive.
Then

\noindent
(i) $\kappa(H_\A)\cap {H''}_K = \kappa(H_k)$

\noindent
(ii) ${H''}_\A = \kappa_\A(H_A)$. 
\end{lem}

\begin{proof}
  Let $x\in \kappa(H_\A)\cap {H''}_K$. Then $\kappa^{-1}(x)$
is a torsor for $H'$. By \cite{harder:75b} this torsor is trivial
which yield the result.

(ii) For each $x\in X$ the map $\kappa_x$, obtained by base change to
$K_x$ is surjective. This is because the Galois cohomology 
$\coh{1}{K_x,H'}$ vanishes since $H'$ is simply connected.
\end{proof}

\begin{pro}
\label{P:compare2}
  Let 
\[
1\rightarrow H' \rightarrow H \stackrel{\kappa}{\rightarrow} {H''} \rightarrow 1
\]
be an exact sequence of connected reductive groups with  $H'$ semisimple simply
connected and  ${H''}$
a torus. Then
\[
\tau(H')\tau({H''}) = \tau(H)
\]
\end{pro}

\begin{proof}
Compare with \cite[Prop 1.2.3]{ono:65}. Let $F$ be compatible
with ${\rm Im}\psi_{H''}$. Consider the integral
\begin{eqnarray*}
  J &=& \int_{H(\A)/H(K)} F(\psi_{H''}(\kappa(x))) d\tau_H \\
    &=& \tau(H) (\sum_{x\in D_q^r} F(x))(\log q)^r.
\end{eqnarray*}
To see this note that $\hat{\kappa}$ induces an isomorphism
on character groups as $H'$ is semisimple. Then apply the above lemma along with
(\ref{P:compatible}). Again by the lemma we may apply
\cite[Theorem 2.4.4]{weil:82} to this integral and obtain:
\begin{eqnarray*}
  J &=& \tau(H')\int_{H(\A)/H(K)} F(\psi_{H''}(y)) d\tau_{H''} \\
    &=& \tau(H')\tau(H'') (\sum_{x\in D_q^r} F(x))(\log q)^r.
\end{eqnarray*}
Again we have made use of the above lemma.
\end{proof}

We now recall Ono's construction of crossed diagrams. Let 
$\tilde{G}$ be the universal cover of $G$, so that we
have an exact sequence
\[
1\rightarrow M\rightarrow \tilde{G} \rightarrow G \rightarrow 1.
\]
Note that $M$ is of multiplicative type, see \cite[IX]{SGA3},
over the field $K$. Recall that

\begin{pro}
  The category of groups of multiplicative type over $K$ 
is antiequivalent to the category of $\gal(\bar{K}^s/K)$-modules
that are finitely generated as abelian groups.
\end{pro}

\begin{proof}
  See \cite[X]{SGA3}.
\end{proof}

The above duality is induced by ${\rm Hom}(-,\gm)$, i.e 
by taking character modules. Now observe that we can find an exact sequence
\[
0\rightarrow M \rightarrow T' \rightarrow T \rightarrow 0,
\]
with $\ch[T'\otimes_K K']$ a projective 
$\gal(K'/K)$-module for some splitting field 
$K'$ of $T'$. To see this set 
\[
G_{M} = \{ g\in\gal(\bar{K}^s/K) | g\text{ fixes }\hat{M} \}.
\] 
Let $\Gamma = \Z[G_{M}]$ and we can find an exact sequence 

\[
0\rightarrow \text{kernel} \rightarrow \Gamma+ \Gamma + \ldots +\Gamma 
\rightarrow M \rightarrow 0.
\]
Set $G^* = (\tilde{G} \times T')/M$ and we have a diagram 
\[
\xymatrix{
   &    & 0  \ar[d]&   &   \\
   &    & T' \ar[d]^i \ar[dr]&   &   \\
 0 \ar[r]& \tilde{G} \ar[r] \ar[dr]& G^* \ar[r]\ar[d]& T \ar[r]& 0 \\
   &    & G  \ar[d]&   &   \\
   &    & 0  &   &  
}
\]

\begin{proof}(of (\ref{T:ono})) We use the above notations.
As in \cite[Lemma 2.1.1]{ono:65} the vertical column above
has a generic section. Assuming Weil's conjecture that 
$\tau(\tilde{G}) = 1$ we obtain
\[ 
\tau(G) = \frac{\tau(T)|{\rm cok}(\hat{i})|}{\tau(T')}
\]
using (\ref{P:compare1}) and (\ref{P:compare2}). The result (\ref{T:ono}), 
follows now from  arithmetic duality theorems and the arguments  in
\cite[pg.99 - 102]{ono:65}. Also note \cite[Corollary 3.3]{oesterle:84}.
\end{proof}

\begin{cor}
\label{C:tsplit}
  Suppose that the group $G$ is split. Then the sequence
\[
\tau_1(G), \tau_2(G), \ldots
\]
is constant. Here $\tau_n(G)$ is the Tamagawa number of the base
change $G\times_k k_n$.
\end{cor}
\begin{proof}
 See the cited work of Ono, in particular \cite[Theorem 2.1.1]{ono:65}
and \cite[Proposition 4.5.1]{ono:63}.
Essentially, under the stated hypothesis,
the Tamagawa number is the cardinality of the fundamental group
which is stable under base change.
\end{proof}

The remainder of this section will be devoted to studying 
how the Tamagawa number changes under base extensions of the
form $k_n/k$ under various hypothesis. We begin by 
recalling the explicit construction of the localization maps
for $\hat{M}$ in our particular setting.

We view the Galois module $\hat{M}$ as a functor on 
field extensions of $K$ in the usual way. Given a diagram 
of fields
\[
\xymatrix{
 K_1' \ar[r] & K_2' \\
 K_1   \ar[r] \ar[u] & K_2 \ar[u]
}
\]
with the vertical maps being Galois extensions we obtain morphisms
\[
\gal(K_2'/K_2) \rightarrow \gal(K_1'/K_1)
\]
and
\[
\hat{M}(K_1') \rightarrow \hat{M}(K_2').
\]
This gives maps 
\[
\coh{i}{K_1'/K_1, \hat{M}} \rightarrow\coh{i}{K_2'/K_2, \hat{M}}.
\]
In particular we have diagrams
\[
\xymatrix{
 \bar{K}^s \ar[r] &  \overline{K_x}^s \\
 K   \ar[r] \ar[u] & K_x \ar[u]
}
\]
for every $x\in X$. This yields the map
\[
\coh{1}{K,\hat{M}}\rightarrow\prod_{x\in X} \coh{1}{K_x,\hat{M}},
\]
whose kernal is $\sha(K,\hat{M})$.

We record here :
\begin{lem}
  Consider the projection $\pi:X_n\rightarrow X$.
 
\noindent
(i) The map $\pi$ is etale.
 
\noindent
(ii) Let $x\in X$ then $\pi^{-1}(x)$ consists of $\text{g.c.d}(n,\deg x)$
points.
                                                                                
(iii) Suppose $\pi(y)=x$. Then $K_{n,y}/K_x$ is a cyclic
Galois extension of degree
\[
\frac{n}{\text{g.c.d}(n,\deg x)}.
\]
Its Galois group is generated by the Frobenius.
\end{lem}

\begin{proof}
  This is well known. See for example,
\cite{rosen:00}.
\end{proof}

We denote by $K_n$ the function field of $X_n$. We assume
from now on that the is a splitting field $L$ for
$G$ that has field of constants $k$ and further 
$k$ contains all roots of unity of order dividing
the fundamental group of $G$. The field
$L_n$ has its obvious meaning.

\begin{lem}
\label{L:tivact}
Let $y\in X_n$ and denote by $\pi$ the
projection $X_n\rightarrow X$.
Under the above hypothesis we have that 
$\hat{M}(K_n)$ (resp. $\hat{M}(K_{n,y})$) is a
trivial $\gal(K_n/K)$-module 
(resp. $\gal(K_{n,y}/K_{\pi(y)})$-module).
\end{lem}

\begin{proof}
  Note that
\[
\hat{M}(K_n) = \hat{M}(\bar{K}^s)^{\gal(\bar{K}^s/K_n)}.
\]
The group $\gal(K_n/K)$ is cyclic and generated by the
Frobenius. Let $F$ be a lift of the Frobenius to $\gal(\bar{K}^s/K)$. 
As the field of constants of $L$ and $K$ are the same, we may assume $F$ fixes
$L$. Now by the assumption on the roots of unity
we have that $F$ acts on $\hat{M}(\bar{K}^s)$ trivially. 
The result follows. The other case is similar.
\end{proof}

\begin{pro}
\label{P:top}
  Suppose that $G$ has a splitting field with
field of constants $k$ and if $k$ contains
all roots of unity dividing the order of $M$
then
\[
|\coh{0}{K_n,\hat{M}}|.
\]
doesn't depend on $n$.
\end{pro}

\begin{proof}
  Follows from the above lemma.
\end{proof}

\begin{lem}
  The natural maps
\[
\coh{i}{K_n/K,\hat{M}} \rightarrow \prod_{y\in X_n}
\coh{i}{K_{n,y}/K_{\pi(y)},\hat{M}}
\]
is injective. 
\end{lem}

\begin{proof}
  By the Riemann hypothesis for function fields we can find a point 
$x\in X$ with $\deg x$ coprime to $n$. If $y$ 
lifts this point we have that
\[
\gal(K_n/K) \cong \gal(K_{n,y}/K_{x}).
\] 
The result follows from (\ref{L:tivact}).
\end{proof}

\begin{theo}
    Suppose that $G$ has a splitting field with
field of constants $k$ and  $k$ contains
all roots of unity dividing the order of $M$.
Then there is a natural isomorphism
\[
\sha(K,\hat{M})\cong \sha(K_n,\hat{M}).
\]
\end{theo}

\begin{proof}
  We have an inflation - restriction sequence inducing the diagram
\[
\xymatrix{
0   \ar[d]                   &  \\
\coh{1}{K_n/K,\hat{M}}\ar[d] &  \\
\coh{1}{K,\hat{M}} \ar[d] \ar[r]^(.4)l & \prod_{x\in X} \coh{1}{K_y,\hat{M}} \ar[d]\\
\coh{1}{K_n,\hat{M}} \ar[d] \ar[r]^(.4){l_n} & \prod_{y\in X_n} \coh{1}{K_y,\hat{M}} \\
\coh{2}{K_n/K,\hat{M}} & 
}
\]

The previous lemma shows that we have an injection
\[
\sha(K,\hat{M})\hookrightarrow \sha(K_n,\hat{M}).
\]
Let $\alpha\in\sha(K_n,\hat{M})$. Also by the
lemma we can lift $\alpha$ to 
$\tilde{\alpha}\in\coh{1}{K,\hat{M}}$. We need to show that
$\tilde{\alpha}$ is in the subgroup $\sha(K,\hat{M})$ modulo the image
of $\coh{1}{K_n/K,\hat{M}}$.  Let 
$l(\tilde{\alpha}) = (\beta_x)_{x\in X}$. Choose for each 
$x\in X$ a lift $\tilde{x}\in X_n$. Since
$\alpha\in\sha(K_n,\hat{M})$ we have that 
\[
\beta_x\in\coh{1}{K_{\tilde{x}}/K_x,\hat{M}}
\]
for every $x$. The hypothesis imply that
$\hat{M}(\widehat{k(\tilde{x})})$ is a trivial as a 
$\gal(\widehat{k(\tilde{x})}/\widehat{k(x)})$-module. So we may think of 
each $\beta_x$ as a homomorphism 
\[
\beta_x:\gal(K_{\tilde{x}}/K_x)
\rightarrow\hat{M}(K_{\tilde{x}}).
\]
We have an inclusion of Abelian groups
\[
\hat{M}(K_n)=\hat{M}(\bar{K}^s)^{\gal(\bar{K}^s/K_n)}\hookrightarrow
\hat{M}(K_{\tilde{x}}).
\]
As each $\beta_x$ comes from $\tilde{\alpha}$ the 
above homomorphisms factor through $\hat{M}(K_n)$.
Define a homomorphism
\[
\alpha_0:\gal(K_n/K) \rightarrow \hat{M}(K_n)
\]
by 
\[
F \mapsto \tilde{\alpha}(F),
\]
where $F$ is the Frobenius and we are choosing a representing
cocycle for $\tilde{\alpha}$ in the above expression.
Now an easy diagram chase shows that
$\tilde{\alpha} - \alpha_0$ is in $\sha(K,\hat{M})$.
\end{proof}

\begin{cor}
\label{C:stable2}
  Under the hypothesis of the theorem and assuming Weil's
conjecture fro the universal cover of $G$ we have
\[
\tau_n(G) = \tau(G).
\]
for every $n$.
\end{cor}

\begin{proof}
  Combine the theorem with (\ref{T:ono}) and (\ref{P:top}).
\end{proof}

\bibliographystyle{alpha}
\bibliography{./../../latestbib/alggrp.bib,./../../latestbib/mybib.bib,./../../latestbib/sga.bib,./../../latestbib/mybib2.bib}

\end{document}